\documentclass[a4paper,reqno]{article} 
\usepackage[T2A]{fontenc}
\usepackage[utf8]{inputenc}
\usepackage{pifont}
\usepackage[dvips]{graphicx}
\usepackage{amscd,amsmath,amssymb,amsthm,latexsym,mathrsfs,pifont}
\usepackage{longtable,cmap,indentfirst,setspace}
\usepackage{ccaption, caption, subcaption}
\usepackage{tikz,tikz-cd,array,hhline}
\usepackage{color,xcolor,hyperref}
\usepackage{verbatim}

\usepackage{geometry}
\hypersetup{
    unicode=true,          
    pdftoolbar=true,        
    pdfmenubar=true,        
    pdffitwindow=false,     
    pdfnewwindow=true,      
    colorlinks=true,        	
    linkcolor={red!50!black},	
    citecolor={blue!50!black},  
    urlcolor={blue!80!black},    
    filecolor=magenta,          
}

\def\C{\mathbb{C}}
\def\P{\mathbb{P}}

\title{Degenerations, transitions and quantum cohomology}
\author{Sergey Galkin}
\begin{document}
\date{April 28, 2015}
\maketitle
Given a singular variety I discuss the relations 
between quantum cohomology of its resolution and smoothing.
In particular, I explain how toric degenerations helps 
with computing Gromov--Witten invariants,
and the role of this story in \href{http://www.fanosearch.net}{``Fanosearch''} programme \cite{fanosearch}.
The challenge is to formulate enumerative symplectic geometry of complex $3$-folds
in a way suitable for extracting invariants under blowups, contractions,
and transitions.

\bigskip

Topologically \emph{a conifold transition} is a surgery on a real $6$-dimensional manifold
switching a $3$-sphere with a $2$-sphere. The boundary of the tubular neighbourhood
of a $2$-sphere or a $3$-sphere is isomorphic to $S^2 \times S^3$,
so one can replace a tubular neighbourhood of one of those spheres with a tubular neighbourhood
of another one. Smith--Thomas--Yau \cite{STY} show that this procedure makes
sense in the context of symplectic manifolds, replacing a Lagrangian $S^3$
with a symplectic $S^2$. There are two natural ways of choosing a $S^2$ related by
a \emph{flop}. Note that similar procedure in real dimension $4$,
namely replacing a Lagrangian $S^2$ with a symplectic $S^2$,
does not change the topology of underlying manifold (as was observed already by Brieskorn),
and can be compensated merely by changing the symplectic form on the same underlying manifold.
For the conifold transition of $6$-dimensional manifolds topology obviously changes,
since the Euler number of two manifolds differs by $2 = e(S^2)-e(S^3)$.

A local algebro-geometric picture for the conifold transition looks as follows.
Consider a function $f : \C^4 \to \C$ given in some coordinates $x,y,z,w$
by a non-degenerate quadratic form
\[ f = (xy - zw) \]
A $3$-fold conifold singularity (or an ordinary double point in dimension $3$)
is locally given by equation $(f=0)$, the fiber of map $f$ over $0$
\[ X_0 := \{ u\in\C^4 : f(u) = 0\} = \{ (x,y,z,w)\in\C^4 : x y = z w \} \]
It has \emph{a smoothing} $X_t$ for $t\in \C$ given as a fiber of $f$ over $t$
\[ X_t := \{ u\in\C^4 : f(u) = t\} = \{ (x,y,z,w)\in\C^4 : xy = zw + t \} \]
and two \emph{small resolutions} 
$\hat{X_0} = \overline{\Gamma_\phi}, \hat{X_0}' = \overline{\Gamma_{\phi'}}$,
given as Zariski closures of the graphs of rational maps $\phi,\phi' : X_0 \to \C\P^1$
\[ \phi (x,y,z,w) := (x:z) = (w:y) \]
\[ \phi' (x,y,z,w) := (x:w) = (z:y) \]
Here all $X_t$ are isomorphic as complex manifolds
\footnote{Thanks to a $\C^*$-action on the total space $(x,y,z,w) \to (\lambda x,\lambda y,\lambda z,\lambda w)$, compatible with the natural action on the base $t \to \lambda^2 t$.},
and as a real $6$-manifold $X_1$ is isomorphic to $T^* S^3$,
the total space of the (co)tangent bundle on $S^3$
\footnote{It is easier to see in coordinates $z_1,\dots,z_4$, where function $f = xy-zw$
has a form $f = \sum_{k=1}^4 z_k^2$ and for real positive $t$.
Then if $z_k = x_k + i y_k$, two vectors in $\mathbb{R}^4$ $x=(x_1,\dots,x_4)$ and $y=(y_1,\dots,y_4)$
are pairwise-orthogonal and square of norms differ by $|x|^2 - |y|^2 = t$}.
The small resolutions $\hat{X_0}$ and $\hat{X_0}'$ are isomorphic as abstract complex manifolds
to the so-called \emph{local $\C\P^1$}, that is the total space of the bundle
$\mathcal{O}_{\C\P^1}(-1)^{\oplus 2}$.

The global picture is similar: we look for a singular complex threefold $Y_{sing}$,
such that $\rm Sing Y_{sing}$ equals to $N$ ordinary double points.
It always has $2^N$ small resolutions which we denote by $Y_{res}$,
however some of them may fail to be quasi-projective.
The question of existence of a smoothing
\footnote{A smoothing of $Y_{sing}$ is a flat projective morphism $f: \mathcal{Y}\to \Delta$
to a disc $\Delta$, such that $f^{-1}(0) = Y_{sing}$ and for some $t\in\Delta$ the fiber
$f^{-1}(t)$ is a smooth complex threefold $Y_{sm}$. In this situation
we say that $Y_{sm}$ degenerates to $Y_{sing}$.}
is more subtle: versal deformation space
for the conifold singularity is smooth and described above, but in some situations
there are local-to-global obstructions. It was observed that for projective manifolds
the existence of a projective resolution and of a smoothing are \emph{mirror dual}
to each other. Friedman \cite{Fr} shows that a smoothing always exists and unique
(that is the versal deformation space is smooth) for Fano threefolds $Y_{sing}$.
For Calabi--Yau threefolds $Y_{sing}$ a smoothing exists iff there is a linear relation
between exceptional symplectic $2$-spheres in $Y_{res}$ \cite{Fr},
and a projective resolution exists iff there is a linear relation between
vanishing Lagrangian $3$-spheres in $Y_{sm}$.
For proper toric varieties $Y_{sing}$ Gelfand--Kapranov--Zelevinsky demonstrate
the existence of a projective small resolution.

Since in complex dimension $2$ conifold transitions do not affect the topology,
they can be used to compute Gromov--Witten invariants and find mirror potentials
for non-toric del Pezzo surfaces of degrees $5$ and $4$. These surfaces
have toric degenerations with $A_1$ singularities, and the crepant resolutions
of these degenerations are smooth toric weak del Pezzo surfaces. There are variety of methods
to compute holomorphic curves and discs on toric surfaces (equivariant method of Givental,
tropical method of Mikhalkin, Cho--Oh, Fukaya--Oh--Ohta--Ono, Chan--Lau--Leung, etc),
and the results of these computations can be transferred back to non-toric del Pezzo.
Also in dimension $2$ it is quite easy to construct symplectic birational invariants,
for example a surface $S$ is rational iff its quantum cohomology $QH(S)$ is
generically semi-simple \cite{Bay}
\footnote{However I do not know any a priori proof that semi-simplicity
persists under the birational contractions, even in dimension $2$.
As an illustration of the problem: non-algebraic $K3$ surface does not
have any non-trivial Gromov--Witten invariant, but
its blowup in a point does --- an invariant counting curves in an exceptional class!}.

For the threefolds situation is quite different. There are theorems relating
quantitative aspects of Gromov--Witten theories of a manifold and its blowup \cite{Gath,Bay}
\footnote{On practice it is often preferable to use less direct ways of computation
as in \cite{CCGK}, thanks to their more compact collection of the enumerative data.}
or conifold transition \cite{LR,NNU,BG,IX}. However none of them is suitable yet
to demonstrate numerous qualitative relations we expect to hold:
birational/transition invariance of symplectic rational connectedness,
non-triviality of GW invariants, sharp birational invariants.
In the remainder I discuss what we know and what we would like to know
about the blowups and transitions of threefolds in the context of ``birational symplectic topology''.

In my thesis \cite{Gal}
I described all conifold transitions from Fano threefolds to toric threefolds,
thus answering a question posed by Batyrev in \cite{Bat}.
The case when $Y_{sing}$ is a toric Fano is arguably the simplest: by Friedman
a smoothing $Y_{sm}$ exists and its topology is unique, so one just have to identify it,
by computing some invariants and using Fano--Iskovskikh--Mori--Mukai's classification.
There are $100$ conifold transitions from smooth Fano to smooth toric weak Fano,
and approximately one half ($44$) of all non-toric Fanos have at least one such transition,
on average two. To a toric $Y_{sing}$ (or $Y_{res}$) one can associate a Laurent polynomial
$W = \sum_v z^v$ 
where the summation runs over all vertices $v$ of the fan polytope of $Y_{sing}$
(equiv. $Y_{res}$, since the resolution is small).
Batyrev put a conjecture that constant terms of powers $W^d$ are equal
to Gromov--Witten invariants $\langle \psi^{d-2} [pt] \rangle_{0,1,d}$,
that count rational curves on $Y_{sm}$ passing through a generic point
with ``tangency conditions'' prescribed by a power of psi-class.
Nishinou--Nohara--Ueda shown that in fact $W$ above is the Floer potential function,
that counts pseudo-holomorphic discs of Maslov index two on $Y_{sm}$
with a boundary on a Lagrangian torus, obtained as a symplectic transport of
a fiber of a moment map $Y_{sing} \to \Delta$. With Bondal in \cite{BG}
we proved a conjecture of Batyrev by expressing the described GW invariant
as a polynomial of the numbers of holomorphic discs, using the relation
between holomorphic curves, passing through a point, and holomorphic discs
bounded on a Lagrangian torus, in the tropical limit.
Same invariants were also computed in \cite{CCGK}
by alternative (and less uniform) methods for all Fano threefolds.
Also Batyrev and Kreuzer described degenerations of Fano threefolds to
nodal half-anticanonical hypersurfaces in toric fourfolds, and it turns out
that almost every Fano has such a degeneration.

In \cite{LR} Li and Ruan give a relation between Gromov--Witten theories of
two threefolds linked by a conifold transition, and very recently Iritani and Xiao
in \cite{IX} reformulated it as a relation between two quantum connections
\footnote{This work prompted me to speak on this subject here.}:
the quantum connection of $Y_{sm}$ (on cohomology of even degree)
is obtained from the quantum connection of $Y_{res}$ as a residue of a sub-quotient.
Similar descriptions could (and should) be also obtained for the quantum connections
of the blowups.

I expect that some sharp birational invariants could be extracted from Gromov--Witten
theory of threefolds \footnote{For more details I refer to my recent talk
\url{http://www.math.stonybrook.edu/Videos/Birational/video.php?f=20150411-1-Galkin}
in Stony Brook on ``New techniques in birational geometry'' conference.},
and that they should be related to the monodromy group of the quantum connection.
Much work still has to be done to understand better the relation between
these monodromies, and to extract an invariant out of it.

Finally, the relation between conifold transitions and its effect on classical topology
should be clarified further. A hypothesis usually referred to as ``Reid's dream'' \cite{Reid}
says that \emph{all} Calabi--Yau threefolds may be connected by a network of conifold transitions.
There is an extensive experimental evidence towards this --- most of known Calabi--Yau threefolds
were shown to be connected by a network. But there are some classical invariants
of the transition, such as the fundamental group
\footnote{Both $S^3$ (resp. $S^2$) has codimension at least $3$ in $Y_{sm}$ (resp. $Y_{res}$),
hence $\pi_1(Y_{sm}) = \pi_1(Y_{sm} \setminus \bigcup S^3_k) =\pi_1(Y_{res} \setminus \bigcup S^2_k)
= \pi_1(Y_{res})$.}. There are a couple of families of non-simply-connected
threefolds without holomorphic $1$-forms and $K=0$ (e.g. $16$ families were
found by Batyrev and Kreuzer as crepant resolutions of hypersurfaces in toric fourfolds),
so these threefolds are clearly disconnected from the conjectural ``simply-connected network''.
Smith--Thomas--Yau \cite{STY} using conifold degenerations observed that possibly
there are lots of non-K\"ahler symplectic Calabi--Yau threefolds, and Fine--Panov--Petrunin
shown that the fundamental group of symplectic Calabi--Yau threefolds can take infinitely
many different values. Non-simply-connected manifold could not be rationally connected,
non-rationally-connected manifold cannot be symplectically rationally connected,
non-symplectically-rationally-connected manifolds cannot have semi-simple quantum cohomology,
which means at least some element in $QH(Y)$ of $Y$ with $\pi_1(Y) \neq 0$ should be nilpotent.
This suggests an interesting question: \emph{given a non-contractible loop $\gamma$ on $Y$,
associate to it a nilpotent in $QH(Y)$}. One possible approach to this could be
via Givental's interpretation \cite{Giv} of the quantum $D$-module of $Y$ with $S^1$-equivariant Floer
theory on a free loop space $\mathcal{L}Y$. We find amusing the contrast
with Seidel's representation $Sei: \pi_1(\rm Symp Y) \to QH(Y)^*$ that associates
an invertible element to a loop in symplectomorphisms \cite{Sei}.

\providecommand{\arxiv}[1]{\href{http://arxiv.org/abs/#1}{arXiv:#1}}
\providecommand{\hdoi}[2]{\href{http://dx.doi.org/#1}{#2}}

\end{document}